\documentclass[11pt,twoside]{amsart}
\usepackage[utf8]{inputenc}
\usepackage[T1]{fontenc}
\usepackage{mathtools}
\usepackage{graphicx}
\usepackage{tikz}
\usetikzlibrary{chains}
\usepackage{tcolorbox}

\usepackage{mathrsfs}

\usetikzlibrary{matrix}
\usetikzlibrary{shapes.multipart}
\usetikzlibrary{arrows}

\PassOptionsToPackage{pdfusetitle,pagebackref,colorlinks}{hyperref}
\usepackage{bookmark}
\hypersetup{
  linkcolor={red!70!black},
  citecolor={green!70!black},
  urlcolor={blue!80!black}
}

\usepackage{amsmath}
\usepackage{amsthm}
\usepackage{amssymb}

\usepackage[all]{xy}

\usepackage{hyperref}

\usepackage[shortlabels]{enumitem}
\usepackage{tikz-cd}

\numberwithin{equation}{section}

\usepackage[OT2,T1]{fontenc}
\DeclareSymbolFont{cyrletters}{OT2}{wncyr}{m}{n}
\DeclareMathSymbol{\Sha}{\mathalpha}{cyrletters}{"58}

\newcommand{\Gr}{{\rm Gr}}

\newcommand{\nilp}{{\rm nilp}}

\newcommand{\QQ}{\mathbb{Q}}

\renewcommand{\to}{\xymatrix@1@=15pt{\ar[r]&}}
\renewcommand{\rightarrow}{\xymatrix@1@=15pt{\ar[r]&}}
\renewcommand{\leftarrow}{\xymatrix@1@=15pt{&\ar[l]}}
\renewcommand{\mapsto}{\xymatrix@1@=15pt{\ar@{|->}[r]&}}
\renewcommand{\twoheadrightarrow}{\xymatrix@1@=18pt{\ar@{->>}[r]&}}
\renewcommand{\hookrightarrow}{\xymatrix@1@=15pt{\ar@{^(->}[r]&}}
\newcommand{\hook}{\xymatrix@1@=15pt{\ar@{^(->}[r]&}}
\newcommand{\congpf}{\xymatrix@L=0.6ex@1@=15pt{\ar[r]^-\sim&}}

\usepackage{graphicx}
\usepackage[capitalise]{cleveref}
\newtheorem{thm}{Theorem}[section]% theorem counter resets every \subsection

\newtheorem{cor}[thm]{Corollary}
\newtheorem{defn}[thm]{Definition}

\newtheorem{conj}[thm]{Conjecture}

\theoremstyle{definition}

\theoremstyle{definition}

\newtheorem{exa}[thm]{Example}
\newtheorem{rmk}[thm]{Remark}
\crefname{thm}{Theorem}{Theorems}
\Crefname{thm}{Theorem}{Theorems}
\Crefname{thm}{Theorem}{Theorems}
\Crefname{thm}{Theorem}{Theorems}
\crefname{lem}{Lemma}{Lemmas}
\Crefname{lem}{Lemma}{Lemmas}
\crefname{Conjecture}{Conjecture}{Conjectures}
\Crefname{Conjecture}{Conjecture}{Conjectures}
\crefname{Corollary}{Corollary}{Corollaries}
\Crefname{Corollary}{Corollary}{Corollaries}
\crefname{Claim}{Claim}{Claims}
\Crefname{Claim}{Claim}{Claims}
\crefname{Proposition}{Proposition}{Propositions}
\Crefname{Proposition}{Proposition}{Propositions}
\crefname{Remark}{Remark}{Remarks}
\Crefname{Remark}{Remark}{Remarks}
\crefname{Definition}{Definition}{Definitions}
\Crefname{Definition}{Definition}{Definitions}
\crefname{Example}{Example}{Examples}
\Crefname{Example}{Example}{Examples}
\crefname{Exercise}{Exercise}{Exercises}
\Crefname{Exercise}{Exercise}{Exercises}

\def\O#1.{\mathcal {O}_{#1}}			
\def\pr #1.{\mathbb P^{#1}}				
\def\af #1.{\mathbb A^{#1}}			
\def\ses#1.#2.#3.{0\to #1\to #2\to #3 \to 0}	
\def\xrar#1.{\xrightarrow{#1}}			
\def\K#1.{K_{#1}}						
\def\bA#1.{\mathbf{A}_{#1}}			
\def\bM#1.{\mathbf{M}_{#1}}
\def\bN#1.{\mathbf{N}_{#1}}
\def\bL#1.{\mathbf{L}_{#1}}				
\def\bB#1.{\mathbf{B}_{#1}}				
\def\bK#1.{\mathbf{K}_{#1}}			
\def\subs#1.{_{#1}}					
\def\sups#1.{^{#1}}						
\def \neone#1.{\overline{{\rm NE}(#1)}}

\begin{document}
\title[Perverse-Hodge octahedron]{Perverse-Hodge octahedron}

\author[M.~Mauri]{Mirko Mauri}
\address[Mirko Mauri]{CNRS UMR7586 Institut de Math\'{e}matiques de Jussieu-Paris Rive Gauche}
\email{mauri@imj-prg.fr}

\begin{abstract}
    The perverse-Hodge octahedron is a 3D enhancement of the Hodge diamond of a compact hyperk\"{a}hler manifold. Its existence is equivalent to Nagai's conjecture, which holds for all known deformation types. The octahedron appears implicitly in \cite{HM} and \cite{ShenYin2023}.
\end{abstract}

\maketitle
{\let\thefootnote\relax\footnotetext{This note was presented at the Séminaire de Géométrie Algébrique at IMJ-Paris Rive Gauche on June 15, 2023 and in occasion of the conference "Hodge theory, tropical geometry and o-minimality" at the Berlin-Brandenburg Academy on October 30, 2023.}}
\marginpar{}

\setcounter{tocdepth}{1}

The sequence of Betti numbers of a compact orientable manifold $X$ is palindromic by Poincar\'{e} duality. If $X$ is further a K\"{a}hler manifold, the sequence of Betti numbers can be rearranged in the shape of a diamond, the so-called Hodge diamond. Summing along the rows of the diamond, we recover the sequence of Betti numbers, but in addition we can visualize the effect of complex conjugation as the symmetry about the vertical axis of the diamond. This was invisible at the level of Betti numbers. 

If $X$ is further hyperk\"{a}hler, the Hodge diamond admits a 3D enhancement which is expected to be an octahedron. As before, summing along the "rows" of the octahedron, we recover the Hodge diamond, and we visualize the effect of the newly discovered P=F symmetry \cite{ShenYin, ShenYin2023} as the invariance of the  3D enhancement under octahedral transformations. Some of these symmetries were invisible at the level of the Hodge diamond.

\begin{exa}
If $X$ is a K3 surface admiting an elliptic fibration $X \to \mathbb{P}^1$, then the four vertices of the octahedron in middle degree represent the symplectic form, its conjugate, a relative ample class for the elliptic fibration, and the pullback of an ample class on $\mathbb{P}^1$.

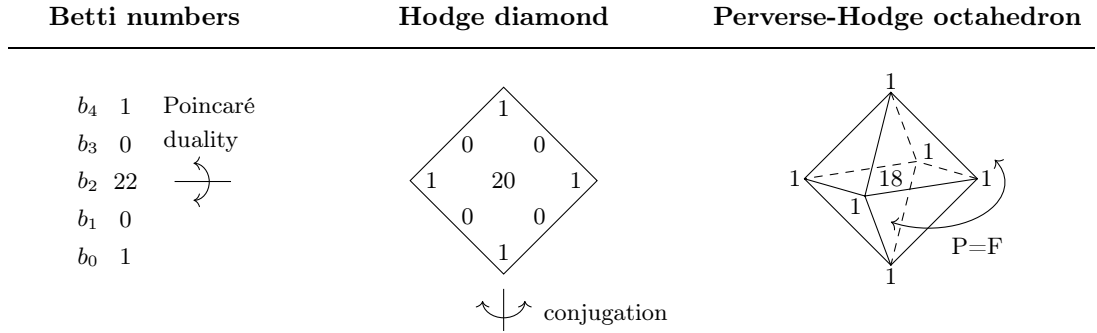
\begin{table}[h]
\centering
\small
\renewcommand{\arraystretch}{2}
   \begin{tabular}{ccccc}
& \textbf{Betti numbers} & $\quad$&\textbf{Hodge diamond} & \textbf{Perverse-Hodge octahedron} \\ \hline \vspace{-1.3 cm} &&\\
& \begin{tikzpicture}[scale = 0.5, every node/.style = {inner sep=1pt, font=\footnotesize}]
\node at (0,0) {$b_0$};
\node at (0,1) {$b_1$};
\node at (0,2) {$b_2$};
\node at (0,3) {$b_3$};
\node at (0,4) {$b_4$};
\node at (1,0) {$1$};
\node at (1,1) {$0$};
\node at (1,2) {$22$};
\node at (1,3) {$0$};
\node at (1,4) {$1$};
%\node at (0,-1) {};
\node at (-1,0) {};
\node[above, text width=1cm, align=center] at (3,2.7) {Poincar\'{e} duality};
\draw[<->] (2.8,1.5) arc(-90:90:12pt and 14pt);
 \draw[] (2.3,2)--(3.8,2);
 \node[above] at (0,-2.3) {};
\end{tikzpicture} & &\begin{tikzpicture}[scale = 0.8, every node/.style = {inner sep=1pt, font=\footnotesize}]
 \node[rectangle,draw,scale=25,rotate=45] at (0,0){}; 
                            \node at (0,-1.2) {1};  
                            \node at (0,+1.2) {1};  
                             \node at (+0.6,+0.6) {0};
                             \node at (+0.6,-0.6) {0};
                             \node at (-0.6,+0.6) {0};
                             \node at (-0.6,-0.6) {0};
                            \node at (0,0) {20};
                            \node at (-1.2,0) { 1};  
                            \node at (+1.2,0) {1};   
                            \node at (-3,0) {};  
                             \node at (+3,0) {};
                             \node at (0,2.25) {};
                             \draw[<->] (-0.35, -2) arc(-180:0:10pt and 8pt);
                             \draw[] (0,-1.8)--(0,-2.5);
                             \node[below] at (1.7,-2) {conjugation};
                             \node[above] at (0,-2.7) {};
                           \end{tikzpicture} & \begin{tikzpicture}[scale=1.15, z={(-.3cm,-.2cm)}, % direction z gets projected to; can also change x,y
                                                                      % use cm to specify non-transformed coordinates
                                                   line join=round, line cap=round, % makes the corners look nicer
                                                   every node/.style = {inner sep=1pt, font=\footnotesize}
                                                  ]
                                 \draw ( 0,1,0) node[above] {1} --
                                       (-1,0,0) node [left] {1} --
                                       (0,-1,0) node[below] {1} --
                                       ( 1,0,0) node[right] {1} --
                                       ( 0,1,0) --
                                       ( 0,0,1) node[below left] {1} --
                                       (0,-1,0) (1,0,0) -- (0,0,1) -- (-1,0,0);
                                 \draw[dashed] (0,1,0) -- (0,0,-1) -- (0,-1,0) (1,0,0) -- (0,0,-1) -- (-1,0,0);
                                 \draw(0,0,0) node[] {18};
                                 \draw(0,0.1,-1.1) node[right] {1};
                                  \node at (0,1.75,0){};
                                  \node[above] at (0,-1.9, 0) {};
                          \draw[ <->] (0, -0.5) arc(-120:30:25pt and 15pt);
                      \node[below] at (1,-0.65) {P=F};
                               \end{tikzpicture}
    \end{tabular}
    \end{table}
\end{exa}

\begin{defn} \cite[\S 1.6]{Del80}
Given a nilpotent endomorphism $N$ of a finite dimensional vector space $V$ of index $\mathrm{nilp}(N) = l$, i.e.\ $N^{l}\neq 0$ and $N^{l+1}= 0$, the \emph{weight filtration of $N$ centered at $l$} is the increasing filtration
$$W_0V\subset W_1V\subset\ldots \subset W_{k}V \coloneqq \sum_{i+j+l=k} \ker N^{i+1} \cap \mathrm{Im} N^{-j} \subset \ldots \subset W_{2l-1}V\subset W_{2l}V=V.$$
%such that $N W_k \subseteq W_{k-2}$, and denoting again by $N$ the induced endomorphism on graded pieces, $N^k \colon \Gr^W_{l+k} V \simeq \Gr^W_{l-k}V$ for every $k\geq 0$, 
%see \cite[\S 1.6]{Del80}. 
\end{defn}

Let $X$ be a compact hyperk\"{a}hler manifold of dimension $2n$. Any cohomology class $\omega \in H^2(X, \mathbb{K})$ defines a nilpotent endomorphism $L_{\omega}$ on $H^*(X, \mathbb{K})$ by cup product. 
\begin{enumerate}
\item 
If $\overline{\sigma} \in H^2(X, \mathcal{O}_X)$ is an anti-symplectic form, then the weight filtration of $L_{\overline{\sigma}}$ on $H^*(X, \mathbb{C})$ is related to the Hodge filtration as follows
\[
W^{\overline{\sigma}}_{k} H^*(X, \mathbb{C}) = \bigoplus_{q \geq 2n-k} H^{p,q}.
\]
\item Consider a Lagrangian fibration $f \colon X \to B$ and let $\beta$ be the pullback of an ample class $\alpha \in H^2(B, \mathbb{Q})$. By \cite[Prop.\ 5.2.4]{deCM2005}, up to renumbering, the weight filtration of $L_{\beta}$ on $H^*(X, \mathbb{Q})$ 
is the so-called perverse filtration
\[W^{\beta}_k H^d(X, \QQ)=P_{d+k-2n}H^{d}(X, \QQ).\]
\end{enumerate}
Since $\beta$ is a class of type $(1,1)$, the graded pieces of the perverse filtration $\Gr^{P}_{k}H^d(X, \mathbb{Q})$ carry a pure Hodge structure of weight $d$ and level at most $d$.

\begin{defn}
The perverse-Hodge octahedron is the 3D table whose entries are the integers
\[
h^{i,k,d} \coloneqq \dim (\Gr^{P}_{d+k}H^{2d}(X, \mathbb{C}))^{d+i, d-i}.
\]
\end{defn}

\begin{rmk}
Although geometrically evocative, the existence of the Lagrangian fibration is not necessary: we can replace $\beta$ with any non-zero isotropic class of type $(1,1)$. %Note in particular that $b_{2}(X) \geq 4$.
\end{rmk}

\subsection{Symmetries of an octahedron}
The (refined) P=F symmetry stands for the equalities $h^{i,k,d}=h^{k,i,d}$. These identities arise from the fact that $\beta$ and $\overline{\sigma}$ are both isotropic classes with respect to the Beauville--Bogomolov quadratic form $q$ on $H^2(X, \mathbb{C})$. With the exception of $b_2=4$ (cf \cite[Rem.\ 2.11]{HM}), the special orthogonal group $SO(H^2(X, \mathbb{C}), q)$ acts transitively on isotropic classes, and lifts to a group of graded algebra automorphisms on the whole cohomology $H^*(X, \mathbb{C})$; see \cite[Prop.\ 2.8]{HM}. In particular, there exists an orthogonal transformation, exchanging $\beta$ and $\overline{\sigma}$, that lifts to an algebra automorphism of $H^*(X, \mathbb{C})$, exchanging the perverse and Hodge filtration and inducing $h^{i,k,d}=h^{k,i,d}$; see \cite[Thm 1.5]{ShenYin2023} and \cite[Cor.\ 3.5]{HM}.

In \cite[Rem.\ 5.3]{ShenYin2023}, Shen and Yin provide a beautiful explanation of why $h^{i, k,d}$ enjoy octahedral symmetries. Let $\omega \in H^2(X, \mathbb{C})$ be a K\"{a}hler class. The Lie algebra generated by the classes of Lefschetz type in $\langle \sigma, \overline{\sigma}, \beta, \omega \rangle$ and the operator $H= \bigoplus (l-2n) \mathrm{id}_{H^l(X)}$ is isomorphic to $\mathfrak{so}(6)$ by \cite{LL, Verb}. The integers $h^{i,k,d} $ are the rank of the eigenspaces for the action of a Cartan subalgebra of $\mathfrak{so}(6)$ on $H^*(X, \mathbb{C})$. The perverse-Hodge octahedron must then be invariant under the action of the Weyl group $S_{4}$ of $\mathfrak{so}(6)$, which acts as the group of rotational octahedral symmetries.

\subsection{Why an octahedron?}
The convex hull of the non-zero $h^{i,k,d}$ (thought as nodes of an integral grid in $\mathbb{R}^3$) lies in a rhombic dodecahedron by \cite[\S 131]{Schnell2023}, but the expectation is that it should actually be an octahedron.

\begin{conj}\label{conj:octahedron}
The convex hull of the non-zero $h^{i,k,d}$ is an octahedron.
\end{conj}
Equivalently, a slice of the convex hull at height $d$ (called perverse-Hodge diamond in \cite[Rem.\ 5.3]{ShenYin2023}) should be the diamond of vertices $(\pm d, 0)$ and $(0, \pm d)$. Note that the slice is certainly inscribed in the square $[-d, d] \times [-d, d]$. Indeed, $h^{i, k, d}=0$ for $|i|>d$ and $|k| > d$, since $\Gr^{P}_{k}H^d(X, \mathbb{Q})$ has level at most $d$ and because of the P=F symmetry. Hence, to prove \cref{conj:octahedron}, it suffices to show that $h^{i, k, d}=0$ for $|i+k|>d$. 
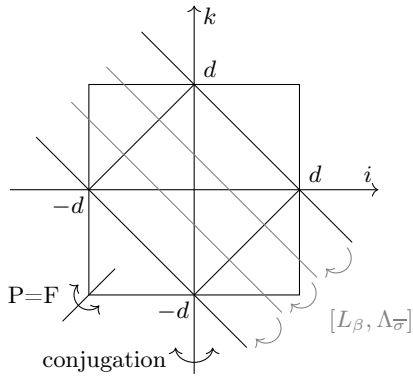
\[
\begin{tikzpicture}[scale = 0.7, every node/.style = {inner sep=1pt, font=\footnotesize}]
\draw[] (0,2)--(-2,0);
\draw[] (0,-2)--(2,0);
\draw[] (2,2)--(-2,2)--(-2,-2)--(2,-2)--(2,2);
\draw[->] (-3.5, 0)--(3.5,0);
\draw[->] (0,-3.5)--(0,3.5);
\node[] at (3.3, 0.3) {$i$};
\node[] at (0.3, 3.3) {$k$};
\draw[] (3, -1)--(-1,3);
\draw[] (-3, 1)--(1,-3);
\draw[gray] (-2.3333, 1.6666)--(1.6666,-2.3333);
\draw[gray] (2.3333, -1.6666)--(-1.6666,2.3333);
                            \node at (0.3,+2.3) {$d$};  
                            \node at (-2.4,-0.3) {$-d$};  
                            \node at (-0.4, -2.3) {$-d$};  
                            \node at (2.3,0.3) {$d$};   
                             \draw[->, gray] (2.9, -1.1) arc(45:-135:8pt and 8pt);
                             \draw[->, gray] (2.9-0.6666, -1.1-0.6666) arc(45:-135:8pt and 8pt);
                             \draw[->, gray] (2.9-0.6666-0.6666, -1.1-0.6666-0.6666) arc(45:-135:8pt and 8pt);
                             \draw[<->] (-2+0.2, -2-0.2) arc(-45:-225:8pt and 8pt);
                             \draw[] (-1.5, -1.5)--(-2.5, -2.5);
                             \draw[<->] (0.35, -3) arc(0:-180:10pt and 8pt);
                             \node[left] at (-2.5, -2) {P=F};
                             \node[right, gray] at (2.5,-2.5) {$ [L_{\beta}, \Lambda_{\overline{\sigma}}]$};
                             \node[below] at (-1.7,-3) {conjugation};
                           \end{tikzpicture}
\]

In \cite[Thm 3.2]{HM}, Huybrechts and the author observe that the sums $\sum_{i+k = \text{const.}}h^{i, k, d}$ are the dimension of the graded pieces of the weight filtration of the nilpotent operator $[L_{\beta}, \Lambda_{\overline{\sigma}}]$, where $\Lambda_{\overline{\sigma}} \coloneqq \lrcorner \sigma$ is the contraction by the symplectic form $\sigma$. Therefore, \cref{conj:octahedron} is equivalent to $\nilp([L_{\beta}, \Lambda_{\overline{\sigma}}]|_{H^{2d}})=d$. This latter statement has a clear geometric counterpart in terms of degenerations of hyperk\"{a}hler manifolds, as explained in the following.

Let $\pi \colon \mathcal{X} \to \Delta$ be a semistable degeneration of hyperk\"{a}hler manifolds deformation equivalent to $X$. For $t \in \Delta^*$, let $N_{2d}$ be the logarithmic monodromy operator on $H^{2d}(\mathcal{X}_{t}, \mathbb{Q})$. Following the general philosophy that the geometry of any compact hyperk\"{a}hler manifold is determined by its Hodge structure of weight two, it was conjectured by Nagai \cite{Nagai2008} that
\[\mathrm{nilp}(N_{2d}) = d \cdot \mathrm{nilp}(N_{2}).\]
\begin{rmk} Recall that $\mathrm{nilp}(N_{2})=\{0,1,2\}$. Nagai's conjecture is open only for $\mathrm{nilp}(N_{2})=1$, proved for all known deformation types of compact hyperk\"{a}hler manifolds, for $d=2n$ and for $2n \leq 8$; see \cite{KLSV, GKLR2022, HM} and references therein. \end{rmk}

\begin{thm}\label{thm:Nagiequivalence} \emph{(\cite[Thm 5.2]{GKLR2022}, \cite[Cor.\ 3.4]{HM})}
If $b_2 \geq 7$, Nagai's conjecture for $\mathrm{nilp}(N_{2})=1$ is equivalent to the existence of the perverse-Hodge octahedron. %; see \cite[Thm 5.2]{GKLR2022} and \cite[Cor.\ 3.4]{HM}.
\end{thm}

Equivalently, Nagai's conjecture would follow from the vanishing 
\[\Gr^P_{i}H^{k,2d-k}(X)=H^{2d-2n}(B, \mathcal{G}_{i,k})=0\]
for $|i-2d+k|>k$, where $\mathcal{G}_{i,k}$ are the coherent sheaves defined in \cite[\S 1.2]{ShenYin2023}; cf also \cite{Schnell2023}.
\cref{thm:Nagiequivalence} traces a new correspondence between degenerations of hyperk\"{a}hler manifolds and Lagrangian fibration, in agreement with the P=W viewpoint \cite{ShenYin, HLSY19:ShenYin}. 
\begin{proof} Let $\pi \colon \mathcal{X} \to \Delta$ be a semistable degeneration of hyperk\"{a}hler manifolds deformation equivalent to $X$ with logarithmic monodromy $N$ and $\mathrm{nilp}(N_{2})=1$ (which exists if $H^2(X, \mathbb{Z})$ contains an isotropic plane, e.g., if $b_2 \geq 7$). The image of $[L_{\beta}, \Lambda_{\overline{\sigma}}]$ and $N$ in $H^2(X, \mathbb{C})$ are isotropic planes. If $b_2 \geq 5$, $SO(H^2(X, \mathbb{C}), q)$ acts transitively on isotropic planes. In particular, an orthogonal transformation exchanges $\mathrm{Im} [L_{\beta}, \Lambda_{\overline{\sigma}}]|_{H^2} = \langle \beta, \overline{\sigma}\rangle$ and $\mathrm{Im} N_2$, and lifts to a graded algebra automorphism which conjugates $[L_{\beta}, \Lambda_{\overline{\sigma}}]$ and $N$; see \cite[Cor.\ 2.10]{HM}. In particular, $\nilp([L_{\beta}, \Lambda_{\overline{\sigma}}]|_{H^{2d}})=\nilp(N_{2d})$. Hence, Nagai's conjecture, i.e., $\nilp(N_{2d})=d$, is equivalent to $\nilp([L_{\beta}, \Lambda_{\overline{\sigma}}]|_{H^{2d}})=d$, equivalently to the existence of the perverse-Hodge octahedron.
\end{proof}

\begin{cor}
All known compact hyperk\"{a}hler manifolds admit a perverse-Hodge octahedron.
\end{cor}

\noindent \textbf{Acknowledgements.} I warmly thank Daniel Huybrechts, Christian Schnell, Junliang Shen and Qizheng Yin for inspiring email exchanges. 
\bibliographystyle{mrl}
\bibliography{HyperHM}
\end{document}